\documentclass[a4paper,fleqn]{cas-sc}

\usepackage[square,numbers]{natbib}

\usepackage[colorinlistoftodos]{todonotes}

\def\tsc#1{\csdef{#1}{\textsc{\lowercase{#1}}\xspace}}
\tsc{WGM}
\tsc{QE}
\tsc{EP}
\tsc{PMS}
\tsc{BEC}
\tsc{DE}

\usepackage{tabularx}
\usepackage{todonotes}

\usepackage{listings}
\usepackage{color}
\definecolor{dkgreen}{rgb}{0,0.6,0}
\definecolor{gray}{rgb}{0.5,0.5,0.5}
\definecolor{mauve}{rgb}{0.58,0,0.82}

\begin{document}
\let\WriteBookmarks\relax
\def\floatpagepagefraction{1}
\def\textpagefraction{.001}
\shorttitle{Portfolio optimisation as ILP}
\shortauthors{C.C.N Kuhn et~al.}

\title [mode = title]{Integer linear programming supporting portfolio design}   

\author[2]{C. C. N. Kuhn}[type=editor,
                        auid=000,bioid=1,
                        orcid=0000-0001-8733-6962]
\ead{carlos.kuhn1@defence.gov.au}
\address[1]{Defence Science and Technology Group, Department of Defence, Canberra, ACT 2610, Australia}
\address[2]{\it{12thLevel Pty Ltd, 14/11 National Circuit, Barton, ACT, 2600}}

\author[1]{G. Calbert}[type=editor,
                        auid=000,bioid=1,
                        ]

\author[1]{I.L Garanovich}[type=editor,
                        auid=000,bioid=1,
                        orcid=0000-0001-6676-6584
                        ]

\author[1]{T. Weir}[type=editor,
                        auid=000,bioid=1,
                        ]

\begin{abstract}
In large organisations and companies, making investment decisions is a complex and challenging task. 
In the Australian Department of Defence (Defence), the complexity is even higher because defence capabilities are public goods and do not have a financial return \textit{per se}. In this work we mathematically define Defence's investment portfolio problem as a Set-Union Knapsack Problem (SUKP). We present a practical way to linearise the model as an Integer Linear Programming (ILP) problem. This linear model was developed as the optimisation engine of the New Investments to Risked Options (NITRO) portfolio selection tool developed by the Defence Science \& Technology Group (DSTG) for Defence force design activities in 2021.
The model is implemented in the Python package called PuLP which can call several linear solver's Application Programming Interface (API), such as GLPK, COIN CLP/CBC, IBM CPLEX, and Gurobi. After comparing the performance of several solvers, we chose Gurobi in the production server. The implementation of the new model and solver enables the rapid execution of exact solutions to the Defence investment portfolio problem. 

\end{abstract}

\begin{keywords}
OR in Defense \sep Integer Linear Programming \sep Portfolio optimisation \sep Combinatorial Optimization
\end{keywords}

\maketitle
\section{Introduction}

The Australian Government has increased its Department of Defence (Defence) funding in recent years and plans to spend close to two percent of Australia's annual Gross Domestic Product (GDP) in the financial year 2022/23. Planned annual spending over the forward estimates exceeds two percent of GDP \cite{ASPI22_23}. In acquiring and maintaining modern and advanced military capabilities, Defence's acquisition and sustainment budgets are planned to double over the coming decade. It is also one of the largest employers in the country. However, changes in the strategic environment, the application of emerging technologies and future uncertainties mean that investment in Defence capabilities has inherent complexities. Defence is a resource intensive enterprise in which expenditure plans are continuously contested. To best deal with uncertainty and complexity, Defence requires an analytical framework for investment planning that allows Defence commanders and leaders to deliver capabilities that mitigate strategic risks in a cost effective manner. 

Defence decision makers are often faced with multiple alternative courses of action, any of which, if implemented, will use different levels of resources and deliver varying benefits. In this work our focus is to help stakeholders make decisions on the Defence investment portfolio in which they have to select a subset of capability proposals with the aim of enhancing Defence's ability to achieve strategic objectives. As one of the Australian Government’s largest departments, Defence manages thousands of projects to support the current force structure and deliver future capabilities. Such projects span the domains of land, air, maritime, space and information and cyber along with analogous enabling domains such as estate and infrastructure, information and communications technology and intelligence. Defence undertakes a continuous review of risks and capability requirements to shape the future direction of the force. Known as the Defence Capability Assessment Program (DCAP), this process is characterised by periodic intense strategic reviews of the force that seek to select a range of capability investments to update its portfolio of projects and mitigate risks.

This problem is shared by practically all organisations seeking to achieve goals by allocating resources in the form of projects and programs.  In the business sector, investments are chosen in order to maximise the company's profit, and the return of the investment is a key measure of benefit. At a personal level, an individual investor seeks to maximise the return of a collection of investments while minimising the risk i.e. the variation in return of those investments.  In Defence, however, we need a different way to measure the investment return because often no agreed-upon measure exists as a decision attribute or because benefits must be couched in relative terms.  Unlike financial instruments, Defence investments are not usually divisible and  generally take a number of years to mature and become effective. But like other investment situations, Defence must assess risk when deciding which investment to purchase and, like the case of a single investor, in order to minimise this risk, Defence consults specialists, in other words subject matter experts (SME). Investors are often supported by financial tools that generate portfolios of investments taking account of risk and expected return.  Similarly, to support the Defence investment problem, a self-consistent mathematical tool to schedule projects that maximise the benefit to Defence while satisfying budget constraints aids decision makers. While the field of portfolio optimisation has been well established over many years, applications in the defence sector deal with the challenges mentioned above and others not seen in other application domains. A thorough survey of portfolio optimisation for defence applications is in \cite{harrison_portfolio_2020}.  This paper demonstrates a novel application of Integer Linear Programming (ILP) to address a particularly challenging portfolio optimisation problem in a defence setting. 




By translating Government and Defence policy guidance into a vision of the future force structure that can be articulated to Government, the DCAP provides a method for Defence to deliver an 'integrated  force by design' \cite{WilliamsFoundation_2017}. The output of the force design process is expressed in the Integrated Investment Program (IIP), which describes Defence's capital investment plan in major equipment, information and communications technology, facilities and workforce. It also programs sustainment funding for in-service and planned investments. The IIP forms the basis for the future force structure. Updates to the IIP to deal with short-term changes and project evolution are presented to Government twice a year for approval in line with the budget cycle. The DCAP has a longer time-horizon and is aimed at developing changes to the IIP in response to strategic and geopolitical changes. The DCAP is an evolving process  and improvements to the DCAP have enabled the force design process to be responsive to changed opportunities and risks \cite{AustralianGovt_2020}. A comprehensive description of the DCAP is in \cite{calbert_2021}.

The DCAP consists of a number of interdependent steps that might vary depending on the portfolio design requirements. Each activity requires SME input and concludes with decisions taken by senior executives. The activities identify changes to the strategic environment, the risks faced by the current and future force structure and develop options to address risks that are tested using various qualitative and quantitative methods. Finally, a selected set of capability options (CO) to address prioritised risks is presented to senior committees for consideration and Government decision.

CO identified in the DCAP consist of a set of projects that collectively address specific risks identified in strategic planning for Defence. The weight, or value of those options, in addressing risks is measured through scenario exploration, war-gaming activities, Bayesian reasoning, SME judgements and simulations. In the DCAP activities in 2021, the principal tool used for CO assessment was through a Bayesian Reasoning Value Model (BRVM) designed to model the current force structure.  The force structure is represented by a hierarchical model consisting of force packages (FPs) containing related lower level capabilities (known as Defence Elements, e.g. a platform such as F-35A Joint Strike Fighter or Infantry Fighting Battalion) which deliver operational effects that contribute to the ability to meet the requirements of  decisive conditions \cite{NATO2010}. Descriptions of the BRVM and CO scoring methodologies are in \cite{nguyen_2021} and \cite{rowe_2021}. The CO represent changes to the force structure and their weights are evaluated through the BRVM as additive multi-attribute utility functions at the decisive conditions at particular points of time within a planning horizon. Other than that high-level description, in this article we are not going to focus on how the weights are measured, but rather how to manage the weighted options. The CO are further grouped into Families, based on the similarities of the effects generated by the CO and the mitigation of risks.  A particular CO may be a scaled version of another. COs are structured so that stakeholders choose to invest in one CO of each family (where each family contains at least the baseline CO of no change to the current portfolio - ie a CO with no projects and value of zero). 

SME make assessments of the value of each option as a function of its effective year, which is the year when all the projects inside the option are ready to contribute to defence goals.  We then need to maximise the sum of the options' values. As we have noticed, we not only choose an option but also determine when the projects in this option are scheduled based on its yearly value. This maximisation is subjected to an annual budget constraint and the cost of an option is determined by the schedule and annual costs of projects within the option.  This means we are not dealing with a simple linear knapsack problem, but with a complex variation of it known as the set-union knapsack problem (SUKP), \cite{Goldschmidt_1994,ARULSELVAN2014214, wu_solving_2020}. 

The knapsack problem is known to be \textit{NP}-hard \cite{CACCETTA20015547} as is its problem variation studied here. A number of algorithms have been developed to solve knapsack problems and their variants.  Popular practical approaches to solve these (either exactly or approximately) include a range of specifically designed methods such as Dynamic Programming and Branch and Bound.   Linear programming relaxation can find approximate solutions. Heuristic approaches such as memetic algorithms have also been studied and are becoming a popular solution strategy \cite{harrison_multi-period_nodate_0}.  Because the knapsack problem is a combinatorial problem, it makes sense to use heuristics to facilitate trade-offs between optimality and computational time. The idea behind the heuristic approach is to randomly generate a population of solutions whose individuals satisfy the constraints of the problem, and then check the value of the portfolio. Combinations of individuals are generated iteratively until a threshold is reached and the algorithm keeps the solution with best value. In using a heuristic, optimality is not guaranteed, i.e., the algorithm is more liked to return a local maximum/minimum. If the optimality is the goal, we need to use other methods that are more computationally demanding.

The SUKP is generalisation of the knapsack problem.  In the SUKP, each item is a set of elements and rather than each element having a value, each item has a nonnegative value.  Each element has a nonnegative weight. A portfolio is a collection of items and the portfolio's weight  is given by the total weight of the elements in the union of the items' sets. The optimal portfolio is one that maximises the values of its included items while remaining within a total weight constraint. This problem has applications to data-base partitioning, machine loading, finance and cryptography. Goldschmidt et al. \cite{Goldschmidt_1994} demonstrates that the SUKP is \textit{NP}-hard, even in very simple cases, develops an exact dynamic programming algorithm and characterises when it will run in polynomial time.  However, for more general problems, no satisfactory approximate solution is usually obtained in polynomial time. More recent papers \cite{ARULSELVAN2014214,wu_solving_2020, Feng_2018,harrison_2021} have focused on heuristic solutions.  Our literature search has not revealed any papers implementing exact ILP solutions to the SUKP.

In this problem, the decision variables are binary, which means  decision variables are represented by 1 or 0 respectively if a CO is in or out of the portfolio. In order to use well established methods from the field of ILP \cite{schrijver_theory_1998}, we propose a novel approach to linearise the SUKP. This means to mathematically define the objective function and the constraints as linear functions of the decision variables. We also investigate ways to relax the problem in case the problem size or its complexity increase. 

This article is presented as follows. We mathematically define the Defence investment portfolio as a set-union knapsack problem at section \ref{model}. We explain its implementation using PuLP, which is a open source python library at section \ref{implementation}. In the section \ref{discussion} we compare some open source solvers with the well known Gurobi \cite{gurobi} commercial solver and discuss the advantages and disadvantages of using ILP for such a task.

\section{Integer Linear Programming Model}\label{model}

In constructing an investment portfolio the following problem structure and features are noted:

(a) A hierarchy of Family-Capability Option-Project is invoked as the basic structure for decision making.  Families consist of multiple COs and each CO consist of one or more projects. Families are aligned with force packages and the COs seek to mitigate risks identified affecting the future force.  Because COs within each family are designed to achieve the same goal, albeit possibly in a different way or possibly to a different level, a maximum of one CO from each family is chosen for entry into the portfolio.

(b) All projects within a selected CO must be included in the portfolio. Projects are indivisible. A project may only be counted once in the portfolio, even though the project may appear in multiple COs and possibly in different families.  

(c) Projects can have multiple variants across different options.  The variation may be in terms of cost spread or duration. Only a single variant of each project may be selected for inclusion in the portfolio.  To avoid the possibility of two variants of the same project being selected, a restriction to this feature is that projects with variants can only appear in one family.

(d) Capability options are valued, rather than individual projects.  This allows for the capture of synergies of projects within options.  The value of each CO is time-dependent and determined by when all projects within the CO can deliver a capability effect.

(e) Each project variant can begin at any year within a predefined time-window of earliest to latest start date. Preferred, earliest and latest start dates are defined for each project.

(f) The value of the whole portfolio is computed as the sum of the values of the options chosen.

(h) The portfolio is constrained by the budget available in each fiscal year. The model allows for under-programming and over-programming in each year.

(g) Dependencies between projects can be explicitly accounted for via the option construct. An option consisting of two projects will necessarily include both projects in the portfolio if the option is selected. Similarly, construction of two options in one family, one consisting of one project and one consisting of two projects, will ensure only one or both projects are included in the portfolio.  

The Defence's investment portfolio problem is defined as an ILP in the following way.
COs are mathematically represented by its binary decision variable $x^o_i$, that means $x^o_i=1$ if the option $i$ is chosen to be in the portfolio or zero otherwise. Each option in the CO set has a scalar value as a function of the latest effective year of projects in the option $i$, $v_{i}(t_{eff})$.
The objective of the optimisation model is to select the subset of COs (and by inference their projects) that together generate the maximum value to Defence. So the problem can be mathematically expressed as:

\begin{equation}
  \max \left( \sum_{i \in CO} v_{i}(t_{eff}) x^o_i \right)
  \label{objective_func}
\end{equation}
subject to:
\begin{equation}
  \sum_{i=1}^{N^o_f} x^o_i \leq 1
  \label{one_opt_contraint}
\end{equation}
The first constraint defined by equation (\ref{one_opt_contraint}) is the \textbf{One option per family constraint}. Notice here the sum is over the options inside of each family $f$ where $N^o_f$ represent the number of options $i$ in each family. This constraint will force the solver to select only one option from each family or that family is not in the solution at all. 

To define the remaining constraints in this ILP problem, we need to understand the content of each option. Each option $i$ will have a group of projects, each of which can be investment or divestment or a combination of those. Each project $j$ has a yearly cost $c_{jt}$ associated to it, a positive cost signalling an investment or negative cost signalling a divestment. Mathematically those projects will be another binary decision variable in the ILP denoted by $x^p_j$. As the decision in eq. (\ref{one_opt_contraint}) is to make in the options, not at project level, those new decision variables will not directly be part of the objective function, eq. (\ref{objective_func}) , but linked to it by constraints. 

The Defence Budget is also a function of the year and is expressed as $b_t$. We then need to make sure the options in the portfolio will satisfy the budget constraint. To do so we add to the ILP the \textbf{Budget Constraint} as

\begin{eqnarray}
   \sum_{j=1}^{N^p} c_{jt} x^p_j\leq\ b_t + \varepsilon \qquad \forall t \in \mathcal{T}
   \label{budget_constraint1}
\end{eqnarray} 
\begin{eqnarray}
   \sum_{j=1}^{N^p} c_{jt} x^p_j\geq\ b_t - \varepsilon  \qquad \forall t \in \mathcal{T}
  \label{budget_constraint2}
\end{eqnarray}
where the sum is over all projects for each year $t$ in the set $\mathcal{T} = \{2022, 2023, ...\}$, where $N^p$ is the total number of projects. We also define a slack, $\varepsilon$, on the the yearly budget, to give an extra flexibility to the constraint. The slack is useful to help the stakeholder to explore whether it would be possible to find extra funding or modify the budget if particularly attractive options are just outside the budget envelope. Another approach is to incorporate the slack into the objective function as previously presented in \cite{calbert_2021}.

The job of the optimiser is to schedule those projects in the DCAP expenditure window by choosing the options in each family that maximise the overall value at (\ref{objective_func}) and satisfy the yearly budget (\ref{budget_constraint1}) and ( \ref{budget_constraint2}). 
The way we deal with the scheduling task is by creating pseudo projects for each project in each of their possible start years. 
For example, lets suppose that P1 can start in the range $\{[2021,2023\}$, then we nominate P1 starting at 2021 and create two new projects P1.1 starting at 2022 and P1.2 starting at 2023, each with the same length and cost profile as P1, so the set of pseudo projects is $\mathcal{PS} = \{P1, P1.1, P1.2\}$;observe that P1 is part of the set $\mathcal{PS}$. The next step then is to create an option for each possible combination of projects' start times. Of course, each project can be scheduled only once, which is easy to control by adding the \textbf{schedule constraint}

\begin{equation}
   \sum_{j=1}^{N^p_{ps}} x^p_j\leq 1 \qquad \forall j \in \mathcal{PS}
  \label{schedule_constraint}
\end{equation}
where the sum is over all pseudo projects $j \in \mathcal{PS}$ . In order to facilitate understanding, we present an example of the entry data for the solver that complements the fictitious family data shown at table \ref{CO_table}.

\begin{table}[h!]
\centering
\begin{tabular}{|c | c | c |c| c |c|}
\hline
\textbf{Family} & \textbf{Capability Option} & \textbf{Projects}  & \textbf{Value}\\
\hline
F1 & F1.0 & None & 0  \\
\hline
F1 & F1.1 & P1, P2 & 0.4 \\
\hline
F2 & F2.0 & None & 0\\
\hline
F2 & F2.1 & P2, P3, P4 & 0.5 \\
\hline
F2 & F2.2 & P5, P6  & 0.3 \\
\hline
\end{tabular}
\caption{Fictitious example of the families, capability options and projects data-set. Notice that projects can be shared between options. The CO F1.0 and F2.0 are the baseline options - options without any new projects. One CO from each Family must be included in the portfolio.}
\label{CO_table}
\end{table}

We need to link projects to the options. In this example we have two families and five options to choose from, and it has six projects. The total number of decision variables will be $N^o + N^P = 5 + 6 = 11$. For example, if option F1.1 is chosen from family F1, we need a constraint that forces the projects P1 and P2 to be also chosen by the algorithm. This we call the \textbf{Project option constraint} and mathematically is defined by

\begin{equation}
   \sum_{j=1}^{N^P_{i}} x^p_j\geq N^P_{i} x^o_i
  \label{proj_opt_constraint}
\end{equation}
where $N^P_i$ is the number of projects in the option $i$. For example, from option F1.1 we have $N^P_i=2$ and then $x^p_1+x^p_2 - 2x^o_2\geq 0$.
When $x^o_2 = 1$ the only way to satisfy the inequality is having $x^p_1 = x^p_2=1$, and when $x^o_2 = 0$ the $x^p_1$ and $x^p_2$  are still undetermined, the projects are free to be used by another option in case we have options sharing same projects, like the options F1.1 and F2.1 in our working example.

This constraint works to link the projects to its options, but we also need to tell the algorithm to choose only projects inside of an option. To do this we add the \textbf{Option project constraint}, where $N^o_{j}$ are the options containing project $j$ 

\begin{equation}
   \sum_{i=1}^{N^o_{j}} x^o_i\geq x^p_j
  \label{opt_proj_constraint}
\end{equation}

Notice, in this case, the sum is over the options $i$ that contain the project $j$. This constraint completes the link in both ways between COs and projects.

Equations (\ref{objective_func}) through (\ref{opt_proj_constraint}) serve to define a portfolio optimisation selection model for our set-union problem in which projects can be shared across options and families.  
Note that a common linear programming specification for this model would invoke the use of additional binary variables $w_k$ and a large number $B$ as constraints to control the requirement that a project $x_{ijk}$ can only be used once as a resource 

\begin{equation}
   \sum_{i=1}^{N^{f}}\sum_{j=1}^{N^{o}} x_{ijk}\leq Bw_k
  \label{add_constrain}
\end{equation}
and the use of an extra term in the objective function {${-}{\sum_{k=1}^{N^p_{k}} w_k}$}  to ensure that the additional binary variables are correctly set when undetermined by constraints; see e.g. \cite{bradley_1977} for specifications and examples to set up the ILP in this way. 

Equations (\ref{objective_func}) through (\ref{opt_proj_constraint}) represent a novel and elegant model for this set-union knapsack optimisation problem that avoids the use of additional binary variables and additional terms in the objective function.  

\begin{table}[ht]
\centering 
\begin{tabularx}{\linewidth}{|l|X|X|X|X|X|X|X|}
\hline
\textbf{Project} & \multicolumn{1}{c|}{\textbf{Mandated}} & \textbf{Fixed in time} & \textbf{Preferred Start Year} & \textbf{Earliest Start Year} & \textbf{Latest Start Year} & \textbf{Year 1} & \textbf{Year 2}
\tabularnewline \hline
P1 & True & False & year 1 & year 1 & year 3 & 3000 & 2300 \tabularnewline 
\hline
P2 & False & True & year 2 & year 1 & year 4 & 0 & 1500 \tabularnewline 
\hline
P3 & False & True & year 1 & year 1 & year 1 & -660 & -1000 \tabularnewline 
\hline
\end{tabularx}
\caption{Fictitious example of projects' data-set format. The Mandated column specifies whether a project must be included in the portfolio. Projects start in the window determined by the Earliest Start Year and Latest Start Year. Each project can also be fixed in time by the user, which means the project will start on the Preferred Start Year if the project is chosen to be in the portfolio. Each project has an annual spend profile for each year shown in columns Year1 and Year2. Projects can be new investments, changes (in start time or expenditure profile) to current investments or divestment of currently scheduled projects. For example P3 is the cessation of a currently scheduled project and therefore has a negative  spend profile as it will return available resources for use on a new project. Divestment projects are by construction fixed in time because they are changes to existing scheduled projects.}
\label{project_table}
\end{table}

Table \ref{project_table} describes three projects for potential inclusion in a portfolio.  The table has two columns with boolean inputs that allow the user to nominate projects that must be included in the portfolio and also projects that are fixed in time. For example project P1 has to be in the portfolio because its Mandated value is True, but it can start at year 1, or year 2, or year 3. Project P2 is not mandated but it is fixed in time, which means, if it is selected, it must start at the Preferred Start year, year 2. 

This optimisation model needs to incorporate the human in the loop, so that it can incorporate portfolio features required by users. We can integrate this portfolio specific settings in a few constraints. 
For example, a mandated project is dealt by adding one extra constraint, the \textbf{Mandated Project constraint} to the ILP problem:

\begin{equation}
   \sum_{j=1}^{N^p_{ps}} x^p_j = 1
  \label{mandated_project}
\end{equation}
This is just an overwrite to the schedule constraint (\ref{schedule_constraint}), which converts the inequality to an equality for the pseudo projects that are mandated. Projects that are fixed in time do not have additional pseudo projects created. 

A simple constraint is also required for mandated options and mandated families. We just set the decision variable of the option to be one in the case of a mandated option. To deal with the mandated family we just overwrite the constraint (\ref{one_opt_contraint}) to become an equality for the mandated decision variable in a manner similar to that made in constraint (\ref{mandated_project}).

As noted in subparagraph (b), this model deals with the value of COs as function of time, which means time dependency also affects the value of the options in table \ref{CO_table}. To calculate option values we also create pseudo options that are the combinations of all pseudo projects. For example, from table \ref{project_table}, project P1 can start at year 1, year 2 or year 3; and P2 can start only at year 2 because it is fixed in time. That means option F1.1, which has projects P1 and P2, will evolve to pseudo-options F1.1.0, F1.1.1 and F1.1.2, where the difference is on the start year of the project P1 to be year 1, year 2 and year 3 respectively. 
The value of those pseudo options will depend on the year that all projects on the option became effective. 
That means that the total number of decision variables to represent a CO will be the product of the number of all possible start years for each project in the CO. In other words, it is the product of the number of pseudo projects for each project in the CO. 
It is noticeable that the ILP problem size increases very fast as the number of projects increases, and for that reason we need to find ways to relax the problem as we discuss at section \ref{discussion}. 

\section{Implementation}\label{implementation}

In this section we explain how we implemented the ILP problem using a open source modeler in python called PuLP \cite{mitchell_optimization_nodate} in the Defence environment. The general process involves the following actions. Templates have been designed in Microsoft Excel for the collection of family, CO and project data and the options for assessment are developed by representatives of the military Services and other departmental and force design experts. A dedicated costing team estimates project costs. The CO are  evaluated and scored by SME who test the options under a specific scenario and time frames through a decision conference, a structured experiment or multiple criteria decision analysis process; this has included applications of the Analytical Hierarchy Process \cite{saaty80}, Best-worst method \cite{rezaei} and other tools including Bayesian network models \cite{nguyen_2021}.  The complete data set is transformed and loaded into a Microsoft SQL Server database. We access the stored data to feed the solver by calling executable stored procedures in the database. The stored procedures return tables with data that is formatted in a similar way to that presented in tables \ref{CO_table} and \ref{project_table}.

\begin{figure}[htbp]
    \caption{NITRO application user interface}
    \includegraphics[width=\textwidth]{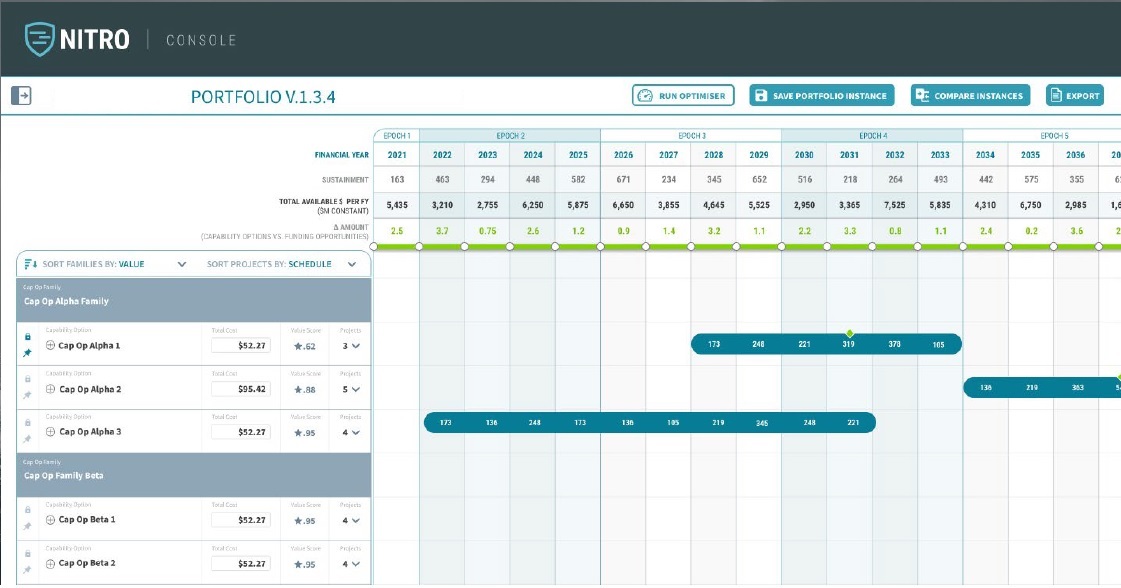}
    \label{Nitro_screen_shot}
\end{figure}

The python script connects to the server using the pyodbc package. The front-end browser user interface (UI) (Figure \ref{Nitro_screen_shot}) developed by 12thLevel Pty Ltd calls the python script to optimise a version of portfolio.  The python code clones the portfolio and returns the optimised result to a table in the Microsoft SQL Server database. The UI displays families, COs and projects in a hierarchical format together with the available budget by year and the cost profile of each project.  Multiple views of the data are available; for example the family view shows the full hierarchy of families, COs and projects and the CO view shows those options that have been selected when the portfolio has been optimised.  The UI allows users to mandate options for inclusion in the portfolio or to disable options so that they are excluded from consideration.  It also allows users to modify project costs and timings or lock a project so that no changes may be made to it. The UI and optimisation process have been designed for the tool to be used in a workshop environment, in which SMEs are able to view portfolio development "on the fly", and investigate and structure portfolios with certain features; for example, users could define a portfolio that preferences maritime capabilities over all other capabilities.  

The optimiser algorithm has two python scripts, the solver script that calls the LP solvers and the python class \textbf{ModelDCAP}. \textbf{ModelDCAP}  contains the module called \textbf{data\_frame\_creator} that cleans, arranges and validates the data to create all the pseudo projects and pseudo options to feed the ILP problem. In this class we also define the ILP problem in the module \textbf{LP\_define}.

The first part of the module \textbf{data\_frame\_creator} simply ensures the raw data is grouped into families with the correct structure like the data presented in table \ref{CO_table}.  

The next part of the module validates the data, making sure all the projects commence on the preferred start date. This procedure also handles the fixed in time projects, setting the Preferred Start Year in the data frame to the year that the user sets it on the graphical user interface (GUI). It also creates all required pseudo projects by shifting the start year of the cost profile of each project within its start year window.
Once we have generated all the projects and pseudo projects by finding all the possible combinations, we create all the options to be added to the data frame and send to the solver.






%
The value of the option is dependent on the time that it becomes effective. The classmethod \textbf{find\_epoch} of the procedure checks the year that projects in an option become effective to determine the option epoch and its value. This method returns the year the option becomes effective and the last project to become effective in the option.

The next important module in the class \textbf{ModelDCAP} is the procedure \textbf{LP\_define} that defines the ILP problem using PuLP. The advantage of using PuLP is that it defines the ILP in a friendly language and can call other solvers using the same ILP file. In this module we implement all the constraints defined in the previous section \ref{model}. As this is the part of the program that sets up the linear program, we explicitly show this python module below:

\begin{lstlisting}
#############################################################################
def LP_define(self, category):
    """Define the ILP problem """
    #Get all the years in the DCAP
    DCAP_years = self.Projects.iloc[0, 10::].index

    ## Start the LP problem definition
    prob = pulp.LpProblem("DCAP_Problem", pulp.LpMaximize)

    #create a list of option
    opt_family = self.FO.Option.tolist()
    opt_project = self.Projects.Project.tolist()
    DCAP_vars = opt_family+opt_project

    #create a dictionary of value for each option and the cost of each project

    value_option = self.option_value.tolist()
    cost_project = self.Projects.iloc[:, 10::].sum(1).tolist()
    #combine the two list to create a dictionary with values for FO and Projects
    total_objective = value_option + cost_project
    #create the dictionary
    value_opt = dict(zip(DCAP_vars, total_objective))
    #define the decision variables
    # category = 'Continuous'
    # category = 'Integer'
    opt_vars = pulp.LpVariable.dicts("Choose", DCAP_vars, 0, 1, cat=category)
    print(f'opt_vars : {len(opt_vars)}')

    # Define the Objective function
    prob += pulp.lpSum([value_opt[i]*opt_vars[i] for i in opt_family])

    #######################################################################
    ################## ADDING THE CONSTRAINTS##############################
    
    #adding to the problem the contraint that each year's budget has some some slack
    #put in matrix form
    cost_matrix = []
    for x in DCAP_years:
        cost_matrix.append(np.array(self.Projects[x].tolist()))

    budget = self.budget.set_index('FinancialYear')
    for ind, row in enumerate(cost_matrix):
        year = int(DCAP_years[ind]))
        if year > budget.index.max():
            break
        prob += pulp.lpSum([row[i]*opt_vars[j] 
            for i, j in enumerate(opt_project)]) <= 
                budget.BudgetValue[year] + budget.MaximumBudgetValue[year] 
        prob += pulp.lpSum([row[i]*opt_vars[j] 
            for i, j in enumerate(opt_project)]) >= 
                budget.BudgetValue[year] + budget.MinimumBudgetValue[year] 

    #######################################################################
    #One project per family constraint and Mandate family constraint
    fo_temp = self.FO.set_index('Family')
    for family in self.Families:
        group = np.array(fo_temp.loc[family].Option)
        if len(group.shape)  == 0:
            continue
        if fo_temp.FamilyMandated.loc[family].unique()[0]:
            prob += pulp.lpSum([opt_vars[j] for j in group]) == 1   
        else:
            prob += pulp.lpSum([opt_vars[j] for j in group]) <= 1      

    ######################################################################
    #Link project to option constraint
    opt_no_base_line = self.FO.Option[self.FO.Projects != 0].tolist()
    proj_NoBaseLine = self.FO.ProjOption[self.FO.Projects != 0].tolist()
    for x, y in zip(proj_NoBaseLine, opt_no_base_line):
        prob += pulp.lpSum([opt_vars[i] for i in x]) >= opt_vars[y]*len(x)

    #####################################################################
    ## Create the constraint to force it to choose only projects within the 
    # family option
    for proj in self.Projects['Project']:
        fo_ind_share_proj = []
        for ind, i in enumerate(self.FO['ProjOption']):
            if i != list(np.zeros(4)) and i != 0:
                for j in i:
                    if proj == j:
                        fo_ind_share_proj.append(ind)

        prob += pulp.lpSum([opt_vars[x] 
            for x in self.FO.Option[fo_ind_share_proj]]) >= opt_vars[proj]

    #####################################################################
    ## Choose the project only once Constraint

    #get the unique ProjectID
    proj_id = np.unique(self.Projects.ProjectID.values)
    proj_temp = self.Projects.set_index('ProjectID')
    #group all project variations (schedule and funding string)
    for proj in proj_id:
        group = np.array(proj_temp.loc[proj].Project)
        if len(group.shape)  == 0:
            continue 
        prob += pulp.lpSum([opt_vars[j] for j in group]) <= 1 

    ########################################################################
    ## Mandate project Constraint
    must_project = self.Projects[self.Projects["Mandated"] == True].copy()
    must_project_group = 
        [must_project.Project[must_project.Project.str.contains(x)].tolist() 
        for x in self.initiatives 
         if must_project.Project[must_project.Project.str.contains(x)].tolist()!= 
                []]

    for proj in must_project_group:
        prob += pulp.lpSum([opt_vars[x] for x in proj]) == 1
        
    ########################################################################
    ## Mandate option Constraint
    must_option = self.FO[self.FO["OptionMandated"] == True].copy()
    original_unique = self.FO.Original_option.unique()
    must_option_group = 
        [must_option.Option[must_option.Option.str.contains(x)].tolist()
        for x in original_unique 
        if must_option.Option[must_option.Option.str.contains(x)].tolist() != []]

    for opt in must_option_group:
        prob += pulp.lpSum([opt_vars[x] for x in opt]) == 1   

    return prob
\end{lstlisting}

\section{Discussion}\label{discussion}
As described in \ref{implementation}, we used PuLP as the modeler to implement the ILP problem and to call the solvers. In this section we compare the performance of the open source COIN-OR Branch-and-Cut (CBC) solver with the state of the art commercial solver Gurobi.

\begin{figure}[htp]
    \centering
    \includegraphics[scale=0.37]{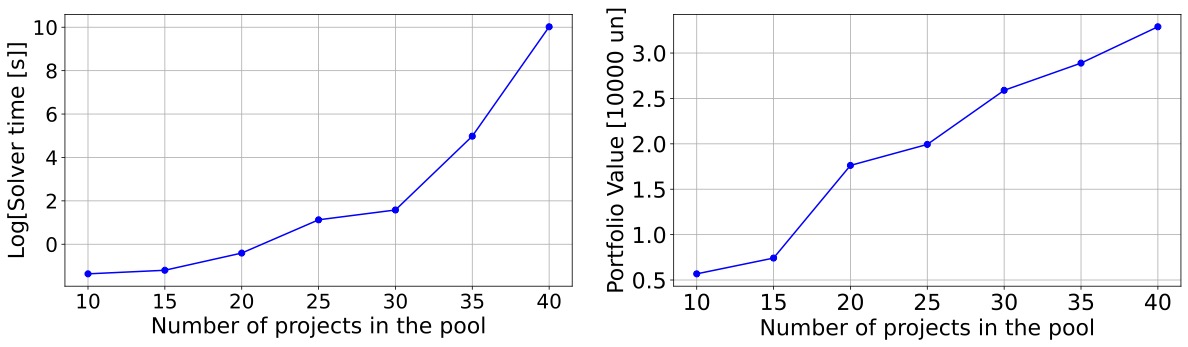}
    \caption{Gurobi solver performance.  The left panel shows the logarithm of the solver time to optimality as a function of number of projects in the pool. The right panel shows the value of the portfolio for the optimal solution (the value of the objective function). The number of decision variables is 16 times the number of projects in the pool, because of the scheduling requirements of the problem.}
    \label{gurobi_int}
\end{figure}

We used simplified simulated data generated using the same procedures created in \cite{harrison_multi-period_nodate_0}. 
Project duration and total cost were independently sampled from a multivariate lognormal distribution found to model defence project expenditures \cite{weir_longitudinal_nodate}.
Here we are assuming that each project has 16 possible start years, and for simplicity each project belongs to one option and each option belongs to one family. That means our test data is a linear knapsack problem. This simplification does not invalidate a comparison of the performance of the solvers as it is a reduction of the SUKP.
It is important to note that the total number of decision variables is 16 times the number of projects in the pool, because of the creation of the pseudo-projects as explained in section \ref{implementation}.
We observed that CBC solver could not solve the ILP problem without any relaxation as it could not find a feasible solution in a reasonable time. 
The Gurobi solver, on the other hand, which uses a mixture of heuristics with Branch-and-Cut algorithms, did solve the same problem and the results are shown in Figure \ref{gurobi_int}. We note that the time to solve increases quickly as function of the problem size.

Our Defence problem deals with a large number of decision variables. Potential portfolio solutions, addressing differing strategies, themes or scenarios, need to be generated quickly in a dynamic workshop environment. To ensure a quick turnaround of portfolios, we investigated ways to find an approximate solution for the problem without compromising the constraints. We found the best way to terminate the optimisation is using the tolerance on the optimality instead of setting a limit run-time.
This means the optimisation will stop when the relative gap between the best known solution and the best known bound on the solution objective is less than a specified tolerance. 
More precisely, if $z_P$ is the primal objective bound (lower bound, best known solution), and $z_D$ is the dual objective bound (the upper bound, normally the solution using continuous variables), then the relative gap is defined as

\begin{equation}
   \mathrm{rel\_gap} = \frac{|z_P - z_D|}{|z_P|}.
  \label{rel_gap}
\end{equation}
This requires a few conditions to be set, and Gurobi  \cite{noauthor_mipgap_nodate} uses the following settings: if $z_P = z_D = 0$, then the gap is defined to be zero. If $z_P = 0$ and $z_D \neq 0$, the gap is defined to be infinity.

In simple terms, the relative gap relaxation compares the feasible known solution with the maximum value possible it could reach if it is allowed to relax the decision variables; in another words, it compares the feasible solution with the continuous relaxed solution. 

One way to conceptualise this process of using continuous relaxation for the upper boundary solution, is to imagine we have a spherical container, and we need to fill the volume with cubes. The smaller the cubes, the more we can fit into the container leaving less empty spaces. The analogy here is the projects are the cubes and the spherical container represents the constraints of the ILP problem with the goal of maximising volume use. 

\begin{figure}[htp]
    \centering
    \includegraphics[scale=0.37]{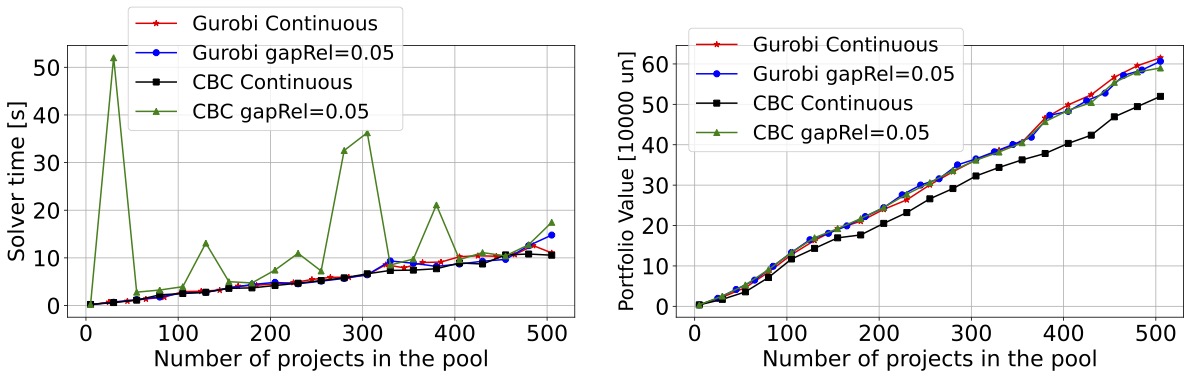}
    \caption{Comparative solver performance. The left panel shows the solver's time as a function of the number of projects in the pool. The right panel shows the value of the portfolio for the optimal solution.}
    \label{gurobi_gap}
\end{figure}

Using those relaxation methods (continuous variable and relative gap), the CBC algorithm and Gurobi present similar performance, where in only a few cases the CBC algorithm takes a much longer time to find a solution.
In the Figure \ref{gurobi_gap} we compare the algorithm performance to find a solution using continuous decision variables $[0,1]$ and relative gap with tolerance set to be 0.05 with integer decision variables $\{0,1\}$. The portfolio value when using continuous variables is calculated by rounding the decision variables, if it is less than 0.5 the variable is set to zero, otherwise it is set to 1, and then summing the value of all option with decision variable equal 1.
The continuous variable solution when using CBC algorithm solver has a lower portfolio value after rounding. This is because the CBC algorithm solution has a higher number of fractional variables increasing the number of continuous variables with value lower then 0.5. It is important to note that we could have slightly higher portfolio values after rounding, but this solution is not feasible, which means it does not satisfy all the constraints of the ILP problem.

Gurobi's superior performance becomes evident when the budget constraint is tighter, the slack defined in equations (\ref{budget_constraint1}) and (\ref{budget_constraint2}) is close to zero, and when we decreased the value of the tolerance.

We also ran both relaxation methods in the SUKP. For this variation of knapsack problem we can not use the continuous variables approach, because the solution might not satisfy the family and options constraints.  A possible approach would be to use mixed integer programming.

The advantage of using relative gap over any continuous variables approach with a lower and upper bound, $[0,1]$, is that the best known solution is feasible, i.e., the solution is not guaranteed to be a global optimum but satisfies all the constraints in the ILP problem.

\section{Conclusion}\label{conclusion}

In this article we have explained the mathematics that define the Defence's Investment Portfolio as a  set-union Knapsack problem. We show that this problem can be expressed as an ILP model and describe how to implement this model using an open-source python linear modeler called PuLP. 
We also compare the performance of the open-source COIN-OR CBC solver with the commercial Gurobi solver. We presented a discussion on relaxation methods to solve the ILP problem. The CBC solver could not solve a simulated portfolio selection problem without some form of relaxation.  This made it unsuitable to use with the general SUKPs that are presented in Defence.
This implementation works as the optimisation engine of the New Investment to Risked Options (NITRO) tool used in Defence portfolio analysis. 

\section{CRediT authorship contribution statement}
\textbf{Carlos Kuhn}: Conceptualisation, Methodology, Software, Data curation, Writing- Original draft preparation.  \textbf{Greg Calbert}: Conceptualisation, Writing- Reviewing and Editing.  \textbf{Ivan Garanovich}: Conceptualisation, Writing- Reviewing and Editing. \textbf{Terence Weir}: Conceptualisation, Data curation, Writing- Revising and Editing, Formal analysis.



\section{Funding}
This work was funded by the Defence Science and Technology Group (DSTG) of the Australian Department of Defence.


\bibliographystyle{elsarticle-num.bst}
\bibliography{references.bib}

\end{document}